\documentclass[hidelinks,onefignum,onetabnum]{siamart251216}

\usepackage{amsmath}
\usepackage{amssymb}
\usepackage{graphicx}
\usepackage{url}

\title{Physics-Informed Neural Networks: A Didactic 
Derivation of the Complete Training Cycle}

\author{
  Abdeladhim Tahimi\thanks{
    Campus de Engenharias e Ci\^{e}ncias Agr\'{a}rias 
    (CECA), Universidade Federal de Alagoas (UFAL), 
    Rio Largo, Alagoas, 57100-000, Brazil 
    (\email{abdeladhim.tahimi@ceca.ufal.br}). 
    Submitted to the editors \today.}
}

\headers{PINNs: Didactic Derivation of the Training Cycle}
{A. Tahimi}

\ifpdf
\hypersetup{
  pdftitle={Physics-Informed Neural Networks: A Didactic Derivation
  of the Complete Training Cycle},
  pdfauthor={Abdeladhim Tahimi}
}
\fi

\usepackage{booktabs}

\begin{document}

\maketitle


\begin{abstract}
This paper is a step-by-step, self-contained guide to the 
complete training cycle of a Physics-Informed Neural Network 
(PINN) --- a topic that existing tutorials and guides 
typically delegate to automatic differentiation libraries 
without exposing the underlying algebra. Using a first-order 
initial value problem with a known analytical solution as a 
running example, we walk through every stage of the process: 
forward propagation of both the network output and its 
temporal derivative, evaluation of a composite loss function 
built from the ODE residual and the initial condition, 
backpropagation of gradients --- with particular attention to 
the product rule that arises in hidden layers --- and a 
gradient descent parameter update. Every calculation is 
presented with explicit, verifiable numerical values using a 
1-3-3-1 multilayer perceptron with two hidden layers and 22 
trainable parameters. From these concrete examples, we derive 
general recursive formulas --- expressed as sensitivity 
propagation relations --- that extend the gradient computation 
to networks of arbitrary depth, and we connect these formulas 
to the automatic differentiation engines used in practice. The 
trained network is then validated against the exact solution, 
achieving a relative $L^2$ error of 
$4.290 \times 10^{-4}$ using only the physics-informed loss, 
without any data from the true solution. A companion 
Jupyter/PyTorch notebook reproduces every manual calculation 
and the full training pipeline, providing mutual validation 
between hand-derived and machine-computed gradients.
\end{abstract}




\section{Introduction}
\label{sec:introduction}

\paragraph{What is a PINN, in one paragraph.}
Physics-Informed Neural Networks 
\cite{raissi2019physics}, building on earlier work 
embedding differential equations into neural network 
training \cite{lagaris1998artificial}, are a framework in 
Scientific Machine Learning (SciML) 
\cite{karniadakis2021physics} that uses a neural network 
to approximate the solution of a problem governed by known 
physical laws. The key idea is to embed these laws --- 
typically expressed as differential equations, but 
potentially including integral constraints, conservation 
principles, or symmetry conditions --- directly into the 
loss function. The network is trained to minimize the 
residual of the governing equations at a set of 
collocation points, while simultaneously satisfying 
prescribed initial or boundary conditions and, when 
available, incorporating experimental or observational 
data. Since their introduction, PINNs have been applied 
in fluid dynamics \cite{raissi2020hidden}, inverse 
problems \cite{karniadakis2021physics}, and a broad range 
of computational science applications 
\cite{cuomo2022scientific}. The methodology continues to 
evolve, with recent work addressing gradient pathologies 
\cite{wang2024causal}, adaptive loss balancing 
\cite{mcclenny2023self}, and broader methodological 
challenges reviewed in \cite{thawon2026physics}. However, 
while the frontier advances rapidly, the foundational 
training mechanics remain the same --- and it is these 
mechanics that we dissect here.

\paragraph{The pedagogical gap.}
Despite the breadth of this literature, a detailed, 
step-by-step account of what a PINN computes during a 
single training step --- from the forward pass through the 
backward pass to the parameter update --- is not available 
in any single resource. Introductory tutorials have been 
developed for mechanics \cite{katsikis2022gentle}, 
astrophysics and plasma physics \cite{baty2024hands}, and 
classroom settings using notebook-based pedagogical 
frameworks \cite{ogueda2024literate}. At the expert level, 
Wang et al.\ \cite{wang2023expert} provide a comprehensive 
guide to training best practices. Standard deep learning 
references \cite{goodfellow2016deep} cover 
backpropagation thoroughly. Yet in all of these, the 
precise sequence of operations by which a physics-informed 
loss and a backward pass produce the gradient that updates 
each weight is either delegated to automatic 
differentiation libraries or presented without exposing 
the underlying algebra --- in particular, the additional 
structure that arises when a differential equation is 
embedded in the loss function receives no explicit 
treatment. A notable exception is Almqvist 
\cite{almqvist2021fundamentals}, who derives the cost 
function gradients symbolically for a 
single-hidden-layer PINN applied to the Reynolds equation. 
That work demonstrates the power of explicit derivation, 
but with one hidden layer the gradient follows a direct 
chain rule --- the recursive structure that emerges when 
parameter perturbations must propagate through multiple 
hidden layers, accumulating product rules, higher order 
derivatives of the activation function, and fan-out 
summations at each stage, does not arise.

\paragraph{What this paper covers.}
The present paper addresses this gap. We take a 
first-order initial value problem with a known analytical 
solution --- $y'(t) + y(t) = 0$, $y(0)=1$ --- and use it 
as a vehicle to expose every detail of the PINN training 
cycle. A small 1-3-3-1 multilayer perceptron with fixed, 
hand-verifiable parameters serves as our working example. 
Specifically, we cover:

\begin{enumerate}
    \item The formulation of a composite loss function 
    comprising an ODE residual term and an initial 
    condition term, weighted by a 
    hyperparameter~$\lambda$.

    \item The simultaneous forward propagation of the 
    network output $\hat{y}(t)$ and its temporal derivative 
    $\hat{y}'(t)$ through the network, layer by layer. 
    This dual propagation is what makes the additional 
    gradient structure of the backward pass visible.

    \item The computation of parametric gradients via 
    backpropagation, with particular attention to the 
    \emph{product rule} that arises in hidden layers 
    because $\hat{y}'$ depends on the network weights 
    through both the activations and their derivatives.

    \item The parameter update via gradient descent.

    \item General recursive formulas that express the 
    parametric gradients for any weight or bias at any 
    layer in terms of two sensitivity quantities --- one 
    for the activation, one for its time derivative --- 
    propagated forward through the network. These formulas 
    extend the concrete examples to networks of arbitrary 
    depth and connect naturally to the automatic 
    differentiation engines used in frameworks such as 
    PyTorch.

    \item The validation of the fully trained network 
    against the known analytical solution 
    $y(t) = e^{-t}$, demonstrating that the 
    physics-informed loss alone --- without any data from 
    the true solution --- is sufficient to recover the 
    correct physical behavior.
\end{enumerate}

Every calculation in this paper is presented with explicit 
numerical values that the reader can (and should) verify 
by hand. To support this, the paper is accompanied by an 
interactive Jupyter/Google Colab Notebook that serves a 
dual purpose: it reproduces every manual calculation using 
PyTorch's automatic differentiation engine, and it 
independently executes the full training and validation 
pipeline, generating the results and figures reported 
here. The exact agreement between the hand-derived values 
and PyTorch's output serves as mutual validation: the 
theory verifies the code, and the code verifies the 
theory.

\paragraph{How this paper is organized.}
\Cref{sec:explanatory_example} introduces the IVP, the 
neural network architecture, and the complete forward pass 
leading to the loss function computation. 
\Cref{sec:backpropagation} details the gradient 
computation via backpropagation --- first through worked 
examples at the output and hidden layers, then through 
general recursive formulas that extend to networks of 
arbitrary depth. \Cref{sec:validation} trains the network 
to convergence and validates it against the analytical 
solution. \Cref{sec:conclusion} presents concluding 
remarks and perspectives for further study.


\section{Explanatory Example of PINN Modeling}
\label{sec:explanatory_example}

\begin{figure}[htbp]
\centering
\includegraphics[width=0.95\textwidth]{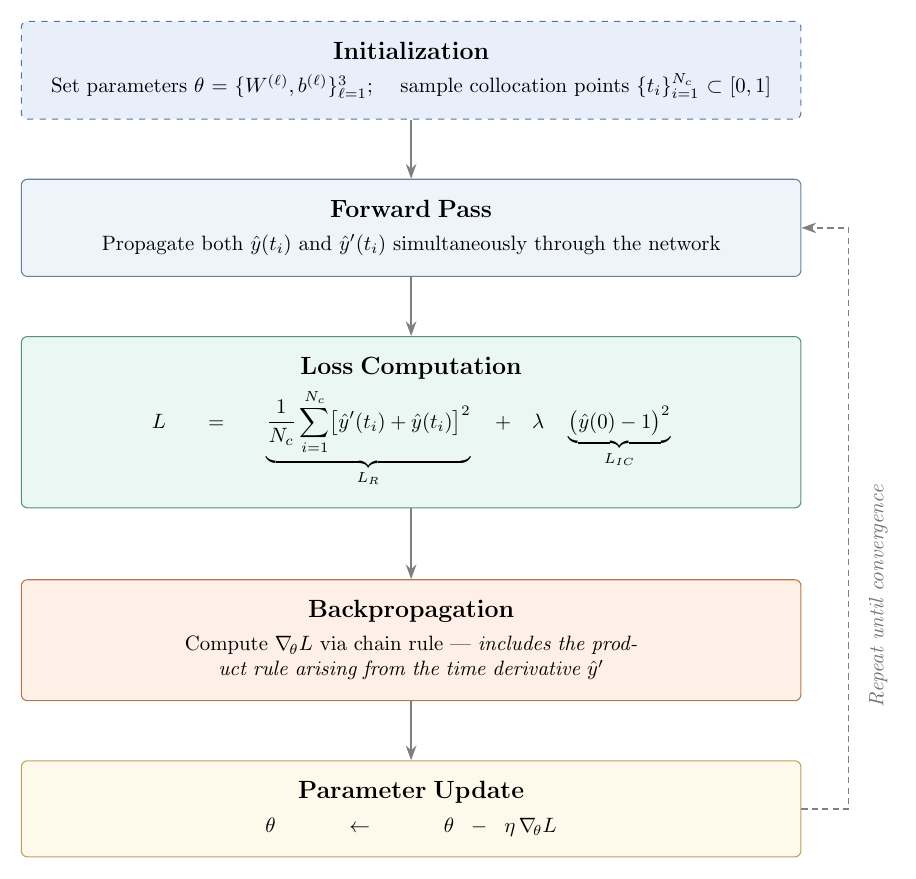}
\caption{The PINN training cycle. After a one-time 
initialization (dashed box), the network iterates 
through four stages: a dual forward pass computing both 
the prediction~$\hat{y}$ and its derivative~$\hat{y}'$, 
evaluation of the composite loss, backpropagation of 
parametric gradients --- where the product rule arising from the 
time derivative $\hat{y}'$ enters in hidden layers --- and a gradient-based 
parameter update. This cycle repeats until the total loss 
converges to a prescribed tolerance.}
\label{fig:training_cycle}
\end{figure}

To demystify the PINN training process, we now work through a 
complete example with every intermediate value written out 
explicitly. \Cref{fig:training_cycle} provides a visual 
overview of the complete training cycle; this section covers 
the forward propagation and loss computation stages, while 
\Cref{sec:backpropagation} details the gradient computation 
and parameter update.

\subsection{Definition of the Initial Value Problem}

Our goal is to find the solution to the first-order ODE
\begin{equation}\label{eq:ode}
y'(t) + y(t) = 0,
\end{equation}
subject to the initial condition $y(0) = 1$. This is the 
``hello world'' of physics-informed learning: simple enough 
to serve as a first exercise with hand-verifiable numbers, 
yet rich enough to expose every mechanism by which the 
network encodes physics --- a residual loss, an initial 
condition constraint, and a derivative that introduces 
the product rule during backpropagation. The analytical 
solution $y(t) = e^{-t}$ will be used in 
\Cref{sec:validation} to assess the accuracy of the trained 
network.

\subsection{Neural Network Configuration}

\begin{figure}[htbp]
\centering
\includegraphics[width=0.9\textwidth]{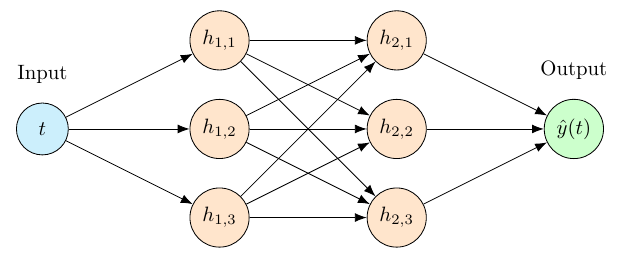}
\caption{Architecture of the 1-3-3-1 MLP. The number of 
trainable parameters for a dense layer is 
$(w_{\mathrm{in}} \times w_{\mathrm{out}}) + w_{\mathrm{out}}$, 
where $w_{\mathrm{in}}$ is the layer input dimension and 
$w_{\mathrm{out}}$ is the number of its neurons. For example, 
the transition between the two hidden layers (both with 3 
neurons) has $(3 \times 3) + 3 = 12$ parameters.}
\label{fig:nn_architecture}
\end{figure}

We employ a \textbf{Multilayer Perceptron (MLP)} with a 
1-3-3-1 architecture: one input neuron (for the variable~$t$), 
two hidden layers with 3 neurons each, and one output neuron 
producing the approximation~$\hat{y}(t)$.
The hidden layers use the hyperbolic tangent ($\tanh$) activation function, and the output layer uses the identity --- a standard choice for regression, since the network 
output must be free to take any real value.
\Cref{fig:nn_architecture} shows the full architecture.

\begin{figure}[htbp]
\centering
\includegraphics[width=0.7\textwidth]{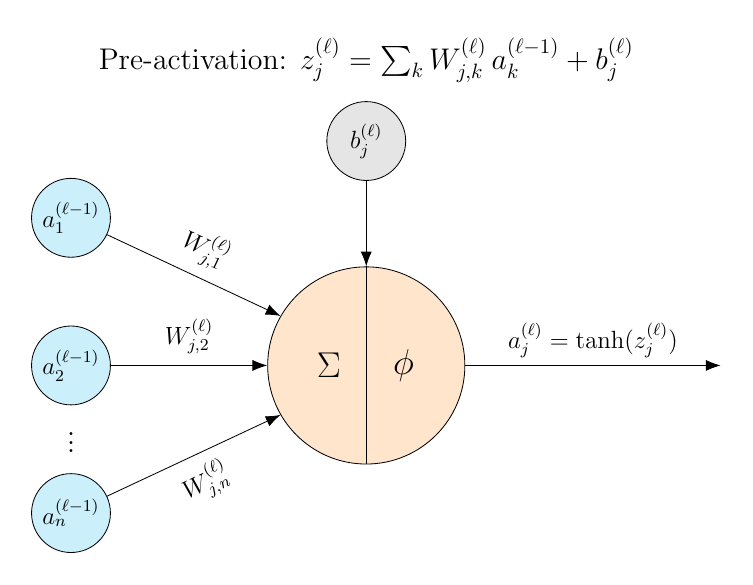}
\caption{Detailed operation of a single artificial neuron. The 
inputs are linearly combined using weights and biases to 
form~$z$, which is then passed through the non-linear 
activation function $\phi(z)$ to produce the output~$a$.}
\label{fig:neuron_detail}
\end{figure}

To fully understand the forward propagation, it helps to look 
at what a single neuron does. As illustrated in 
\Cref{fig:neuron_detail}, a neuron in hidden layer~$\ell$ 
receives outputs from the previous layer, $a^{(\ell-1)}$. It 
performs an affine transformation: multiply the inputs by 
weights~$W^{(\ell)}$, sum the results, and add a 
bias~$b^{(\ell)}$ to produce the pre-activation~$z^{(\ell)}$. 
A non-linear activation function $\phi(\cdot)$ is then applied 
to produce the neuron's output, 
$a^{(\ell)} = \phi(z^{(\ell)})$.

The choice of $\phi(\cdot)$ deserves a careful remark, because 
it is not arbitrary. Without a non-linear activation, every 
layer would compute an affine transformation, and composing 
affine transformations yields another affine transformation --- 
regardless of the number of layers. The network would be 
limited to representing straight lines, and could never capture 
the curvature of a solution like $e^{-t}$. The non-linear 
activation is precisely what breaks this limitation: the 
Universal Approximation Theorem \cite{hornik1989multilayer} 
guarantees that a feedforward network with at least one hidden 
layer and a non-constant, bounded, continuous activation 
function can approximate any continuous function on a compact 
set to arbitrary accuracy, provided sufficiently many neurons. 
Without the non-linearity, this result fails entirely.

For PINNs, however, non-linearity alone is not enough. The 
loss function involves physical derivatives --- $\hat{y}'$ in 
our case, and potentially $\hat{y}''$ or higher for other 
equations --- computed via automatic differentiation through 
the network. This imposes a second requirement: the activation 
function must be \emph{sufficiently smooth}. The Rectified 
Linear Unit (ReLU), the default activation in standard deep 
learning, is non-linear and satisfies the universal 
approximation property, but it is not smooth: it is 
non-differentiable at $z=0$ and its second derivative is zero 
everywhere else. A ReLU-based network simply cannot represent 
the physics of second-order differential equations. PINNs 
therefore require activations that are both non-linear and 
continuously differentiable. The output layer is the 
exception: it uses the identity function $\phi(z) = z$, 
because the network prediction must be free to take any 
real value --- $\tanh$, being bounded to $(-1,1)$, would 
impose an artificial constraint on the solution range.

We use $\tanh$ here, which is the standard activation in the 
PINN literature since the foundational work of Raissi et al.\ 
\cite{raissi2019physics}. It is infinitely differentiable, 
zero-centered --- meaning its output range $(-1,1)$ is 
symmetric around the origin, which promotes balanced gradient 
flow during training --- and its derivative, 
$\phi'(z) = 1 - \tanh^2(z)$, is both smooth and easy to 
compute analytically. For the smooth, non-periodic solution 
treated here, $\tanh$ is a natural and sufficient choice. 
Other smooth activations have been explored in the PINN 
literature, including sinusoidal activations 
\cite{sitzmann2020implicit}, which are well suited for 
problems with oscillatory or periodic solutions, and 
adaptive activations such as Swish and its variants, which 
have shown improvements in convergence for certain problem 
classes \cite{jagtap2020adaptive}. The choice of activation 
remains an active area of research 
\cite{thawon2026physics}; the key requirement for any PINN 
application is that $\phi(\cdot)$ be smooth enough to support 
the order of differentiation demanded by the governing 
equation.

The network parameters (weights $W$ and biases $b$) are fixed 
for this example to allow numerical verification. In a real 
application, these would be set using a standard initialization 
scheme --- such as Xavier \cite{glorot2010understanding} or 
Kaiming \cite{he2015delving} --- which draw random values from 
distributions carefully designed to maintain stable gradient 
flow across layers. The training configuration in 
\Cref{sec:validation} uses Kaiming uniform initialization. The values 
are given in 
\Cref{tab:params_l1,tab:params_l2,tab:params_l3}.

\begin{table}[htbp]
    \centering
    \footnotesize
    \caption{Parameters of Hidden Layer 1.}
    \label{tab:params_l1}
    \begin{tabular}{@{}cc@{}}
        \toprule
        Weights $W^{(1)}$ & Biases $b^{(1)}$ \\
        \midrule
        $\begin{pmatrix} 0.2 \\ -0.5 \\ 0.8 \end{pmatrix}$ &
        $\begin{pmatrix} -0.1 \\ 0.3 \\ -0.2 \end{pmatrix}$ \\
        \bottomrule
    \end{tabular}
\end{table}

\begin{table}[htbp]
    \centering
    \footnotesize
    \caption{Parameters of Hidden Layer 2.}
    \label{tab:params_l2}
    \begin{tabular}{@{}cc@{}}
        \toprule
        Weights $W^{(2)}$ & Biases $b^{(2)}$ \\
        \midrule
        $\begin{pmatrix}
             0.1 & -0.3 &  0.5 \\
             0.6 &  0.2 & -0.4 \\
            -0.2 &  0.7 &  0.1
        \end{pmatrix}$ &
        $\begin{pmatrix} 0.2 \\ -0.1 \\ 0.4 \end{pmatrix}$ \\
        \bottomrule
    \end{tabular}
\end{table}

\begin{table}[htbp]
    \centering
    \footnotesize
    \caption{Parameters of the Output Layer (Layer 3).}
    \label{tab:params_l3}
    \begin{tabular}{@{}cc@{}}
        \toprule
        Weights $W^{(3)}$ & Biases $b^{(3)}$ \\
        \midrule
        $\begin{pmatrix} 0.9 & -0.6 & 0.3 \end{pmatrix}$ &
        $\begin{pmatrix} -0.3 \end{pmatrix}$ \\
        \bottomrule
    \end{tabular}
\end{table}

\subsection{Loss Function Computation: Forward Propagation}

The total loss function $L$ is what guides the training, and it 
has two parts:

\begin{enumerate}
    \item \textbf{Residual Loss ($L_R$):} This term enforces 
    the ODE. It is computed over a set of collocation points 
    $\{t_i\}_{i=1}^{N_c}$. At each point~$t_i$, the residual 
    is $R(t_i) = \hat{y}'(t_i) + \hat{y}(t_i)$ --- which 
    would be zero if the network output satisfied the ODE 
    exactly. The loss is the mean squared residual:
    \begin{equation}\label{eq:loss_R}
    L_R = \frac{1}{N_c} \sum_{i=1}^{N_c} R(t_i)^2.
    \end{equation}

    \item \textbf{Initial Condition Loss ($L_{IC}$):} This 
    term enforces $y(0) = 1$. It is simply the squared error 
    at the initial point:
    \begin{equation}\label{eq:loss_IC}
    L_{IC} = \bigl(\hat{y}(0) - 1\bigr)^2.
    \end{equation}
\end{enumerate}

The total loss is the weighted sum 
$L = L_R + \lambda \cdot L_{IC}$, where $\lambda$ is a 
hyperparameter that balances the two terms. Since our purpose 
here is to illustrate every detail of the training process, we 
keep things as simple as possible: a single collocation point 
$t_c = 0.5$ and $\lambda = 10$. The following subsections trace 
the computation of $\hat{y}$ and $\hat{y}'$ at this point and 
at $t = 0$.

\subsubsection{Forward Propagation for the Residual 
               Loss ($L_R$)}

The first step is to compute the network prediction and its 
derivative at the collocation point $t_c = 0.5$. The network 
input is $a^{(0)} = 0.5$ and --- since the input is simply~$t$ 
--- its derivative with respect to~$t$ is 
$\dot{a}^{(0)} = 1.0$.

\paragraph{Step 1: First Hidden Layer ($\ell=1$).}
We compute the pre-activation $z^{(1)}$ and its temporal 
derivative $\dot{z}^{(1)}$:
\begin{equation}\label{eq:z1_R}
z^{(1)} = W^{(1)} a^{(0)} + b^{(1)}
    \quad \text{and} \quad
    \dot{z}^{(1)} = W^{(1)} \dot{a}^{(0)}.
\end{equation}
The $\tanh$ activation is then applied. A useful identity is 
that the derivative of $\tanh$ can be written in terms of the 
activation itself: $\phi'(z) = 1 - \tanh^2(z)$, so 
$\phi'(z^{(1)}_j) = 1 - (a^{(1)}_j)^2$. The time derivative 
of the activation then follows from the chain rule: 
$\dot{a}^{(1)}_j = \phi'(z^{(1)}_j)\,\dot{z}^{(1)}_j$. 
In vector form:
\begin{equation}\label{eq:a1_R}
a^{(1)} = \tanh(z^{(1)})
    \quad \text{and} \quad
    \dot{a}^{(1)} = \bigl(1 - (a^{(1)})^2\bigr) 
    \odot \dot{z}^{(1)},
\end{equation}
where $\odot$ denotes the Hadamard (element-wise) product. 
Take a moment to verify a few entries in 
\Cref{tab:results_l1_R} before moving on.

\begin{table}[htbp]
\centering
\footnotesize
\caption{Values at the First Hidden Layer for $L_R$ 
($t = 0.5$).}
\label{tab:results_l1_R}
\begin{tabular}{@{}ccccc@{}}
\toprule
Neuron ($j$) & $z_j^{(1)}$ & $\dot{z}_j^{(1)}$
             & $a_j^{(1)}$ & $\dot{a}_j^{(1)}$ \\
\midrule
1 & $0.0000$ & $0.2000$  & $0.0000$ & $0.2000$  \\
2 & $0.0500$ & $-0.5000$ & $0.0500$ & $-0.4988$ \\
3 & $0.2000$ & $0.8000$  & $0.1974$ & $0.7688$  \\
\bottomrule
\end{tabular}
\end{table}

\paragraph{Step 2: Second Hidden Layer ($\ell=2$).}
The process repeats, using the outputs of the previous layer 
as inputs:
\begin{equation}\label{eq:z2_R}
z^{(2)} = W^{(2)} a^{(1)} + b^{(2)},
    \quad
    \dot{z}^{(2)} = W^{(2)} \dot{a}^{(1)},
\end{equation}
\begin{equation}\label{eq:a2_R}
a^{(2)} = \tanh(z^{(2)}),
    \quad
    \dot{a}^{(2)} = \bigl(1 - (a^{(2)})^2\bigr) 
    \odot \dot{z}^{(2)}.
\end{equation}
The results are in \Cref{tab:results_l2_R}.

\begin{table}[htbp]
\centering
\footnotesize
\caption{Values at the Second Hidden Layer for $L_R$ 
($t = 0.5$).}
\label{tab:results_l2_R}
\begin{tabular}{@{}ccccc@{}}
\toprule
Neuron ($k$) & $z_k^{(2)}$ & $\dot{z}_k^{(2)}$
             & $a_k^{(2)}$ & $\dot{a}_k^{(2)}$ \\
\midrule
1 & $0.2837$  & $0.5540$  & $0.2763$  & $0.5117$  \\
2 & $-0.1690$ & $-0.2873$ & $-0.1673$ & $-0.2792$ \\
3 & $0.4547$  & $-0.3123$ & $0.4257$  & $-0.2559$ \\
\bottomrule
\end{tabular}
\end{table}

\paragraph{Step 3: Output Layer ($\ell=3$).}
The output layer uses the identity activation 
($\phi(z) = z$), so the network output $\hat{y}$ equals the 
final pre-activation $z^{(3)}$, and the output derivative 
$\hat{y}'$ equals $\dot{z}^{(3)}$:
\begin{equation}\label{eq:z3_R}
z^{(3)} = W^{(3)} a^{(2)} + b^{(3)}
    \quad \text{and} \quad
    \dot{z}^{(3)} = W^{(3)} \dot{a}^{(2)}.
\end{equation}
\Cref{tab:results_l3_R} gives the numerical values.

\begin{table}[htbp]
\centering
\footnotesize
\caption{Values at the Output Layer for $L_R$ ($t = 0.5$).}
\label{tab:results_l3_R}
\begin{tabular}{@{}cccc@{}}
\toprule
$z_1^{(3)}$ & $\dot{z}_1^{(3)}$ & $\hat{y}$ & $\hat{y}'$ \\
\midrule
$0.1768$ & $0.5513$ & $0.1768$ & $0.5513$ \\
\bottomrule
\end{tabular}
\end{table}

\paragraph{Residual Loss Computation.}
With these values in hand, the ODE residual at $t = 0.5$ is:
\begin{equation}\label{eq:residual}
R(0.5) = \hat{y}'(0.5) + \hat{y}(0.5) 
       = 0.5513 + 0.1768 = 0.7281.
\end{equation}
The residual loss is the square of this value:
\begin{equation}\label{eq:LR_value}
L_R = R(0.5)^2 = (0.7281)^2 \approx 0.5301.
\end{equation}

\subsubsection{Forward Propagation for the Initial Condition 
               Loss ($L_{IC}$)}

Next, we compute the network prediction at $t = 0$ to evaluate 
the initial condition loss. The network input is 
$a^{(0)} = 0$ and its derivative is $\dot{a}^{(0)} = 1$.

\paragraph{Step 1: First Hidden Layer ($\ell=1$).}
Since $a^{(0)} = 0$, the pre-activation reduces to the bias 
vector, and the derivative of the pre-activation equals the 
weight vector:
\begin{equation}\label{eq:z1_IC}
z^{(1)} = W^{(1)} \cdot 0 + b^{(1)} = b^{(1)}
    \quad \text{and} \quad
    \dot{z}^{(1)} = W^{(1)} \cdot 1 = W^{(1)}.
\end{equation}
\Cref{tab:results_l1_IC} collects the results.

\begin{table}[htbp]
\centering
\footnotesize
\caption{Values at the First Hidden Layer for $L_{IC}$ 
($t = 0$).}
\label{tab:results_l1_IC}
\begin{tabular}{@{}ccccc@{}}
\toprule
Neuron ($j$) & $z_j^{(1)}$ & $\dot{z}_j^{(1)}$
             & $a_j^{(1)}$ & $\dot{a}_j^{(1)}$ \\
\midrule
1 & $-0.1000$ & $0.2000$  & $-0.0997$ & $0.1980$  \\
2 & $0.3000$  & $-0.5000$ & $0.2913$  & $-0.4576$ \\
3 & $-0.2000$ & $0.8000$  & $-0.1974$ & $0.7688$  \\
\bottomrule
\end{tabular}
\end{table}

\paragraph{Step 2: Second Hidden Layer ($\ell=2$).}
\Cref{tab:results_l2_IC} presents the results.

\begin{table}[htbp]
\centering
\footnotesize
\caption{Values at the Second Hidden Layer for $L_{IC}$ 
($t = 0$).}
\label{tab:results_l2_IC}
\begin{tabular}{@{}ccccc@{}}
\toprule
Neuron ($k$) & $z_k^{(2)}$ & $\dot{z}_k^{(2)}$
             & $a_k^{(2)}$ & $\dot{a}_k^{(2)}$ \\
\midrule
1 & $0.0039$  & $0.5415$  & $0.0039$  & $0.5415$  \\
2 & $-0.0226$ & $-0.2802$ & $-0.0226$ & $-0.2801$ \\
3 & $0.6041$  & $-0.2830$ & $0.5398$  & $-0.2005$ \\
\bottomrule
\end{tabular}
\end{table}

\paragraph{Step 3: Output Layer ($\ell=3$).}
\Cref{tab:results_l3_IC} gives the final values.

\begin{table}[htbp]
\centering
\footnotesize
\caption{Values at the Output Layer for $L_{IC}$ ($t = 0$).}
\label{tab:results_l3_IC}
\begin{tabular}{@{}cccc@{}}
\toprule
$z_1^{(3)}$ & $\dot{z}_1^{(3)}$ & $\hat{y}$ & $\hat{y}'$ \\
\midrule
$-0.1210$ & $0.5953$ & $-0.1210$ & $0.5953$ \\
\bottomrule
\end{tabular}
\end{table}

\paragraph{Initial Condition Loss Computation.}
Since $\hat{y}(0) = -0.1210$, we get:
\begin{equation}\label{eq:LIC_value}
L_{IC} = \bigl(\hat{y}(0) - 1\bigr)^2
       = (-0.1210 - 1)^2
       = (-1.1210)^2
       \approx 1.2566.
\end{equation}

\subsubsection{Total Loss and the Start of Optimization}

With the forward pass complete, we can now put the pieces 
together.

\paragraph{Total Loss Computation.}
From the calculations above, $L_R \approx 0.5301$ and 
$L_{IC} \approx 1.2566$. With $\lambda = 10$, the total loss 
is:
\begin{equation}\label{eq:total_loss}
L = L_R + \lambda \cdot L_{IC}
  \approx 0.5301 + 10 \cdot 1.2566 = 13.0961.
\end{equation}
This single number tells us how far the network is from a 
correct solution at its current parameter state.

\paragraph{Why we need the weighting.}
In practice, when $L_R$ is evaluated over hundreds or thousands 
of collocation points, it tends to dominate the total loss, 
making $L_{IC}$ negligible. The hyperparameter~$\lambda$ 
corrects this imbalance by amplifying the initial condition 
contribution. Typical values reported in the literature are 10, 
50, 100, 200, and 1000 \cite{wang2021understanding}, and the 
choice is critical for model performance.

\paragraph{Optimization via Gradient Descent.}
The training objective is to find the parameters that 
minimize~$L$. The \textbf{trainable parameters}, collectively 
denoted by 
$\theta = \{W^{(\ell)}, b^{(\ell)}\}_{\ell=1}^{L}$, are 
updated iteratively using \textbf{Gradient Descent}:
\begin{equation}\label{eq:gd_update}
\theta_{\text{new}} = \theta_{\text{old}} 
                    - \eta \nabla_{\theta} L,
\end{equation}
where $\eta$ is the \textbf{learning rate}. We use 
$\eta = 0.01$ here --- small enough to ensure stable 
convergence, and a common default for gradient descent.

\paragraph{The need for backpropagation.}
To apply \cref{eq:gd_update}, we need the gradient 
$\nabla_{\theta} L$. It is worth pausing to note an important 
distinction: the \textbf{temporal derivative} ($d/dt$) computed 
during the forward pass tells us how the output changes with 
the input~$t$; the \textbf{parametric derivative} 
($\partial / \partial \theta$) computed during backpropagation 
tells us how the loss changes with the weights and biases. 
These are fundamentally different quantities, and it is 
essential not to confuse them. The next section details the 
backpropagation computation for our example.


\section{Gradient Computation: Backpropagation}
\label{sec:backpropagation}

With the forward pass complete and the total loss computed, we 
now turn to the backward pass. The goal is to compute the 
gradient of the loss with respect to every trainable parameter, 
$\nabla_{\theta} L$. By linearity of differentiation, the 
total gradient is the weighted sum of the gradients of each 
loss component:
\begin{equation}\label{eq:grad_total}
\nabla_{\theta} L = \nabla_{\theta} L_R +
                    \lambda \cdot \nabla_{\theta} L_{IC}.
\end{equation}
The backpropagation algorithm applies the chain rule 
systematically, moving from the output layer back toward the 
input. We demonstrate the computation for one weight and one 
bias from both the output and hidden layers --- enough to 
expose all the key mechanisms. The same logic applies to every 
other parameter.

\subsection{Gradients for the Output Layer}

\paragraph{Weight $W^{(3)}_{1,2}$.}
This is the weight connecting the second neuron of the 
second hidden layer to the output neuron 
(see \Cref{fig:nn_architecture}). The loss 
$L_R = R^2 = (\hat{y} + \hat{y}')^2$ depends on 
$W^{(3)}_{1,2}$ through both $\hat{y}$ and $\hat{y}'$. By the 
chain rule:
\begin{equation}\label{eq:dLR_dW3}
\frac{\partial L_R}{\partial W^{(3)}_{1,2}} = 2R
    \left( \frac{\partial \hat{y}}{\partial W^{(3)}_{1,2}} +
           \frac{\partial \hat{y}'}{\partial W^{(3)}_{1,2}} 
    \right).
\end{equation}
Since the output layer uses the identity activation, we have 
$\frac{\partial \hat{y}}{\partial W^{(3)}_{1,2}} = a^{(2)}_2$ 
and 
$\frac{\partial \hat{y}'}{\partial W^{(3)}_{1,2}} = 
\dot{a}^{(2)}_2$. For the initial condition loss:
\begin{equation}\label{eq:dLIC_dW3}
\frac{\partial L_{IC}}{\partial W^{(3)}_{1,2}} =
    2\bigl(\hat{y}(0) - 1\bigr)\, a^{(2)}_2.
\end{equation}
\Cref{tab:grad_W3} shows the numerical results.

\begin{table}[htbp]
\centering
\footnotesize
\caption{Gradient computation for weight $W^{(3)}_{1,2}$.}
\label{tab:grad_W3}
\begin{tabular}{lll}
\toprule
\textbf{Component} & \textbf{Formula} & \textbf{Result} \\
\midrule
$\partial L_R / \partial W^{(3)}_{1,2}$ &
$2R \cdot (a^{(2)}_2 + \dot{a}^{(2)}_2)$ &
$2(0.7281)(-0.1673 - 0.2792) = \mathbf{-0.6502}$ \\[4pt]
$\partial L_{IC} / \partial W^{(3)}_{1,2}$ &
$2(\hat{y}-1) \cdot a^{(2)}_2$ &
$2(-0.1210 - 1)(-0.0226) = \mathbf{0.0507}$ \\
\bottomrule
\end{tabular}
\end{table}

\noindent The total gradient is:
\begin{equation}\label{eq:grad_W3_total}
\frac{\partial L}{\partial W^{(3)}_{1,2}} =
    -0.6502 + 10 \cdot (0.0507) = \mathbf{-0.1432}.
\end{equation}

\paragraph{Bias $b^{(3)}_{1}$.}
This is the bias of the output neuron. The computation is analogous, with the simplification 
$\frac{\partial \hat{y}}{\partial b^{(3)}_{1}} = 1$ and 
$\frac{\partial \hat{y}'}{\partial b^{(3)}_{1}} = 0$, since 
the bias is not multiplied by any $t$-dependent variable. 
\Cref{tab:grad_b3} shows the details.

\begin{table}[htbp]
\centering
\footnotesize
\caption{Gradient computation for bias $b^{(3)}_{1}$.}
\label{tab:grad_b3}
\begin{tabular}{lll}
\toprule
\textbf{Component} & \textbf{Formula} & \textbf{Result} \\
\midrule
$\partial L_R / \partial b^{(3)}_{1}$ &
$2R \cdot (1 + 0)$ &
$2(0.7281) \cdot 1 = \mathbf{1.4562}$ \\[4pt]
$\partial L_{IC} / \partial b^{(3)}_{1}$ &
$2(\hat{y}-1) \cdot 1$ &
$2(-1.1210) = \mathbf{-2.2420}$ \\
\bottomrule
\end{tabular}
\end{table}

\noindent The total gradient is:
\begin{equation}\label{eq:grad_b3_total}
\frac{\partial L}{\partial b^{(3)}_{1}} =
    1.4562 + 10 \cdot (-2.2420) = \mathbf{-20.9638}.
\end{equation}

\subsection{Gradients for the Hidden Layers}

Here is where things get interesting --- and where the 
time derivative in the loss introduces additional 
complexity in the gradient computation.

\paragraph{Weight $W^{(2)}_{1,2}$.}
This is the weight connecting the second neuron of the 
first hidden layer to the first neuron of the second 
hidden layer (see \Cref{fig:nn_architecture}). For 
parameters in hidden layers, the influence on the output 
is indirect. The dependence of $\hat{y}$ on 
$W^{(2)}_{1,2}$ follows the path 
$W^{(2)}_{1,2} \rightarrow z^{(2)}_1 \rightarrow 
a^{(2)}_1 \rightarrow z^{(3)}_1 \rightarrow 
a^{(3)}_1 = \hat{y}$, and the chain rule gives:
\begin{equation}\label{eq:dyhat_dW2}
\frac{\partial \hat{y}}{\partial W^{(2)}_{1,2}} =
    \underbrace{\frac{\partial \hat{y}}{\partial z^{(3)}_1}}_{=\,1}
    \cdot
    \underbrace{\frac{\partial z^{(3)}_1}{\partial a^{(2)}_1}}_{=\,W^{(3)}_{1,1}}
    \cdot
    \underbrace{\frac{\partial a^{(2)}_1}{\partial z^{(2)}_1}}_{=\,\phi'(z^{(2)}_1)}
    \cdot
    \underbrace{\frac{\partial z^{(2)}_1}{\partial W^{(2)}_{1,2}}}_{=\,a^{(1)}_2}
    = W^{(3)}_{1,1} \cdot \phi'(z^{(2)}_1) \cdot a^{(1)}_2.
\end{equation}
Now here is the key point. The residual loss depends on 
both $\hat{y}$ and $\hat{y}'$, so the gradient 
$\partial L_R / \partial W^{(2)}_{1,2}$ requires both 
$\partial \hat{y} / \partial W^{(2)}_{1,2}$ --- which we 
just computed in \cref{eq:dyhat_dW2} --- and 
$\partial \hat{y}' / \partial W^{(2)}_{1,2}$. To see 
where this second term leads, recall from the forward 
pass (\cref{eq:z3_R}) that 
$\hat{y}' = W^{(3)}_{1,1}\,\dot{a}^{(2)}_1 + 
W^{(3)}_{1,2}\,\dot{a}^{(2)}_2 + 
W^{(3)}_{1,3}\,\dot{a}^{(2)}_3$. 
Since $W^{(2)}_{1,2}$ only affects the first neuron of 
the second hidden layer, only the $\dot{a}^{(2)}_1$ term 
contributes, and we need 
$\partial \dot{a}^{(2)}_1 / \partial W^{(2)}_{1,2}$. 
From the forward pass (\cref{eq:a2_R}), this activation 
derivative is the product 
$\dot{a}^{(2)}_1 = 
\phi'(z^{(2)}_1)\,\dot{z}^{(2)}_1$. 
The parameter $W^{(2)}_{1,2}$ influences \emph{both} 
factors through two distinct paths:
\begin{itemize}
    \item $W^{(2)}_{1,2} \rightarrow z^{(2)}_1 
    \rightarrow \phi'(z^{(2)}_1)$,
    \item $W^{(2)}_{1,2} \rightarrow \dot{z}^{(2)}_1$.
\end{itemize}
\noindent Since both factors depend on the same parameter, 
differentiating their product requires the product rule:
\begin{equation}\label{eq:product_rule}
\frac{\partial \dot{a}^{(2)}_1}{\partial W^{(2)}_{1,2}} =
    \underbrace{\frac{\partial \phi'(z^{(2)}_1)}{\partial W^{(2)}_{1,2}}}_{=\,\phi''(z^{(2)}_1)\,a^{(1)}_2}
    \cdot\, \dot{z}^{(2)}_1
    \;+\;
    \phi'(z^{(2)}_1) \cdot
    \underbrace{\frac{\partial \dot{z}^{(2)}_1}{\partial W^{(2)}_{1,2}}}_{=\,\dot{a}^{(1)}_2},
\end{equation}
where $\phi''(z) = -2\tanh(z)(1 - \tanh^2(z))$. Note that 
$\phi''$ appears here even though the governing ODE involves 
only a first-order time derivative: the second derivative of 
the activation arises not from the physics, but from 
differentiating the gradient $\phi'(z^{(2)}_1)$ with respect 
to the network parameters during backpropagation. This is a 
concrete illustration of why PINNs require smooth activations 
--- even a first-order equation demands $\phi''$ in the 
gradient computation, and higher-order equations would require 
correspondingly higher derivatives. The full derivative of $\hat{y}'$ with respect to $W^{(2)}_{1,2}$ is therefore:
\begin{equation}\label{eq:dyprime_dW2}
\frac{\partial \hat{y}'}{\partial W^{(2)}_{1,2}} =
    \underbrace{\frac{\partial \hat{y}'}{\partial \dot{a}^{(2)}_1}}_{=\,W^{(3)}_{1,1}}
    \cdot
    \frac{\partial \dot{a}^{(2)}_1}{\partial W^{(2)}_{1,2}}
    = W^{(3)}_{1,1} \left(
        \phi''(z^{(2)}_1)\, a^{(1)}_2\, \dot{z}^{(2)}_1 +
        \phi'(z^{(2)}_1)\, \dot{a}^{(1)}_2
    \right).
\end{equation}
\Cref{tab:grad_W2} presents the numerical evaluation.
\begin{table}[htbp]
\centering
\footnotesize
\caption{Gradient computation for weight $W^{(2)}_{1,2}$.}
\label{tab:grad_W2}
\begin{tabular}{lll}
\toprule
\textbf{Component} & \textbf{Formula} & \textbf{Result} \\
\midrule
$\partial L_R / \partial W^{(2)}_{1,2}$ &
$2R\bigl(\partial\hat{y}/\partial W^{(2)}_{1,2} + 
\partial\hat{y}'/\partial W^{(2)}_{1,2}\bigr)$ &
$1.4562 \cdot (0.0416 - 0.4274) = \mathbf{-0.5619}$ \\[4pt]
$\partial L_{IC} / \partial W^{(2)}_{1,2}$ &
$2(\hat{y}-1)\,\partial\hat{y}/\partial 
W^{(2)}_{1,2}$ &
$-2.2420 \cdot (0.9 \cdot 0.9999 \cdot 0.2913) = 
\mathbf{-0.5878}$ \\
\bottomrule
\end{tabular}
\end{table}

\noindent The total gradient is:
\begin{equation}\label{eq:grad_W2_total}
\frac{\partial L}{\partial W^{(2)}_{1,2}} =
    -0.5619 + 10 \cdot (-0.5878) = \mathbf{-6.4399}.
\end{equation}

\paragraph{Bias $b^{(2)}_{1}$.}
This is the bias of the first neuron of the second hidden 
layer. The derivation follows the same logic as for $W^{(2)}_{1,2}$, 
with the simplifications 
$\frac{\partial z^{(2)}_1}{\partial b^{(2)}_{1}} = 1$ and 
$\dot{z}^{(2)}_1$ independent of $b^{(2)}_{1}$:
\begin{equation}\label{eq:dyhat_db2}
\frac{\partial \hat{y}}{\partial b^{(2)}_{1}} =
    W^{(3)}_{1,1} \cdot \phi'(z^{(2)}_1)
    \quad \text{and} \quad
    \frac{\partial \hat{y}'}{\partial b^{(2)}_{1}} =
    W^{(3)}_{1,1} \cdot \phi''(z^{(2)}_1)\, \dot{z}^{(2)}_1.
\end{equation}
\Cref{tab:grad_b2} shows the results.

\begin{table}[htbp]
\centering
\footnotesize
\caption{Gradient computation for bias $b^{(2)}_{1}$.}
\label{tab:grad_b2}
\begin{tabular}{lll}
\toprule
\textbf{Component} & \textbf{Formula} & \textbf{Result} \\
\midrule
$\partial L_R / \partial b^{(2)}_{1}$ &
$2R\bigl(\partial\hat{y}/\partial b^{(2)}_{1} + 
\partial\hat{y}'/\partial b^{(2)}_{1}\bigr)$ &
$1.4562 \cdot (0.8314 - 0.2546) = \mathbf{0.8400}$ \\[4pt]
$\partial L_{IC} / \partial b^{(2)}_{1}$ &
$2(\hat{y}-1)\,\partial\hat{y}/\partial b^{(2)}_{1}$ &
$-2.2420 \cdot (0.9 \cdot 0.9999) = \mathbf{-2.0176}$ \\
\bottomrule
\end{tabular}
\end{table}

\noindent The total gradient for $b^{(2)}_{1}$ is:
\begin{equation}\label{eq:grad_b2_total}
\frac{\partial L}{\partial b^{(2)}_{1}} =
    0.8400 + 10 \cdot (-2.0176) = \mathbf{-19.3360}.
\end{equation}
%
\subsection{Parameter Update}

The backpropagation algorithm computes these gradients for 
\emph{all} 22 parameters of the network. An optimizer then 
uses them to execute the update step. With learning rate 
$\eta = 0.01$, the update for $W^{(3)}_{1,2}$ is:
\begin{equation}\label{eq:update_example}
W^{(3)}_{1,2,\text{new}} = W^{(3)}_{1,2,\text{old}}
    - \eta \frac{\partial L}{\partial W^{(3)}_{1,2}}
    = -0.6 - (0.01)(-0.1432) \approx -0.5986.
\end{equation}
This cycle --- forward pass (loss computation), backward pass 
(gradient computation), and parameter update --- is repeated 
thousands of times. Each iteration nudges the network 
parameters a small step in the direction that reduces the loss, 
progressively adjusting the network to find a solution that 
satisfies both the ODE and its initial condition.

We demonstrated the gradient computation for specific 
parameters in the output and second hidden layers. But what 
about the remaining parameters --- in particular, those in 
the first hidden layer, where the influence must propagate 
through two subsequent layers? The next subsection presents 
the general formulas that cover any parameter at any layer.

\subsection{General Gradient Formulas}
\label{subsec:general_gradients}

With the full training cycle now complete for specific 
parameters, we step back to see the general pattern. We 
begin with a first hidden layer parameter --- where the 
gradient is more involved than anything we have seen so 
far --- and then extract the general formulas that apply 
to any parameter at any layer.

\subsubsection{A First Hidden Layer Example: Weight 
$W^{(1)}_{2,1}$}

This is the weight connecting the input to the second 
neuron of the first hidden layer 
(see \Cref{fig:nn_architecture}). Its influence on the 
output must propagate through \emph{two} subsequent 
hidden layers before reaching $\hat{y}$ and $\hat{y}'$. 
This is what makes it more complex than the second hidden 
layer parameters treated earlier.

\paragraph{Derivative of $\hat{y}$.}
The dependence of $\hat{y}$ on $W^{(1)}_{2,1}$ starts at 
neuron 2 of layer 1, but then fans out to \emph{all three} 
neurons of layer 2, since $a^{(1)}_2$ feeds into every 
neuron of the next layer. The chain rule therefore 
produces a sum over $k = 1, 2, 3$:
\begin{equation}\label{eq:dyhat_dW1}
\frac{\partial \hat{y}}{\partial W^{(1)}_{2,1}} = 
\sum_{k=1}^{3} 
    \underbrace{\frac{\partial \hat{y}}{\partial z^{(3)}_1}}%
    _{=\,1} \cdot
    \underbrace{\frac{\partial z^{(3)}_1}{\partial a^{(2)}_k}}%
    _{=\,W^{(3)}_{1,k}} \cdot
    \underbrace{\frac{\partial a^{(2)}_k}{\partial z^{(2)}_k}}%
    _{=\,\phi'(z^{(2)}_k)} \cdot
    \underbrace{\frac{\partial z^{(2)}_k}{\partial a^{(1)}_2}}%
    _{=\,W^{(2)}_{k,2}} \cdot
    \underbrace{\frac{\partial a^{(1)}_2}{\partial z^{(1)}_2}}%
    _{=\,\phi'(z^{(1)}_2)} \cdot
    \underbrace{\frac{\partial z^{(1)}_2}{\partial W^{(1)}_{2,1}}}%
    _{=\,a^{(0)}}.
\end{equation}
Compare this with the single-path chain rule for 
$W^{(2)}_{1,2}$ in \cref{eq:dyhat_dW2}: the structure is 
the same, but now there are three paths (one per neuron 
in layer 2), and we must sum over all of them.

\paragraph{Derivative of $\hat{y}'$.}
This is where the complexity grows. Recall that at each 
hidden layer, the time derivative of the activation is 
$\dot{a}^{(\ell)}_i = \phi'(z^{(\ell)}_i)\,
\dot{z}^{(\ell)}_i$, and we saw in \cref{eq:product_rule} 
that differentiating this product with respect to a 
parameter that affects both factors requires the product 
rule. For $W^{(1)}_{2,1}$, this product rule applies at 
\emph{two} layers --- both layer 1 and layer 2 --- because 
the parameter's influence propagates through both.

\paragraph{Product rule at layer 1.}
At layer 1, the situation is identical to what we saw for 
$W^{(2)}_{1,2}$ at layer 2. The parameter $W^{(1)}_{2,1}$ 
influences both factors of 
$\dot{a}^{(1)}_2 = \phi'(z^{(1)}_2)\,\dot{z}^{(1)}_2$ 
through two paths:
\begin{itemize}
    \item $W^{(1)}_{2,1} \rightarrow z^{(1)}_2 
    \rightarrow \phi'(z^{(1)}_2)$,
    \item $W^{(1)}_{2,1} \rightarrow \dot{z}^{(1)}_2$.
\end{itemize}
\noindent Applying the product rule:
\begin{equation}\label{eq:pr_layer1}
\frac{\partial \dot{a}^{(1)}_2}{\partial W^{(1)}_{2,1}} =
    \underbrace{\frac{\partial \phi'(z^{(1)}_2)}%
    {\partial W^{(1)}_{2,1}}}%
    _{=\,\phi''(z^{(1)}_2)\,a^{(0)}}
    \cdot\, \dot{z}^{(1)}_2
    \;+\;
    \phi'(z^{(1)}_2) \cdot
    \underbrace{\frac{\partial \dot{z}^{(1)}_2}%
    {\partial W^{(1)}_{2,1}}}%
    _{=\,\dot{a}^{(0)}}.
\end{equation}

\paragraph{Product rule at layer 2.}
Here is the new element. At layer 2, the product 
$\dot{a}^{(2)}_k = \phi'(z^{(2)}_k)\,\dot{z}^{(2)}_k$ 
also depends on $W^{(1)}_{2,1}$ --- indirectly, through 
the layer 1 outputs. Both factors are again affected:
\begin{itemize}
    \item $W^{(1)}_{2,1} \rightarrow a^{(1)}_2 
    \rightarrow z^{(2)}_k \rightarrow \phi'(z^{(2)}_k)$,
    \item $W^{(1)}_{2,1} \rightarrow \dot{a}^{(1)}_2 
    \rightarrow \dot{z}^{(2)}_k$.
\end{itemize}
\noindent Applying the product rule at this layer:
\begin{equation}\label{eq:pr_layer2}
\frac{\partial \dot{a}^{(2)}_k}{\partial W^{(1)}_{2,1}} =
    \phi''(z^{(2)}_k)\,
    \underbrace{\frac{\partial z^{(2)}_k}%
    {\partial W^{(1)}_{2,1}}}%
    _{=\,W^{(2)}_{k,2}\,P^{(1)}_2}
    \cdot\, \dot{z}^{(2)}_k
    \;+\;
    \phi'(z^{(2)}_k) \cdot
    \underbrace{\frac{\partial \dot{z}^{(2)}_k}%
    {\partial W^{(1)}_{2,1}}}%
    _{=\,W^{(2)}_{k,2}\,Q^{(1)}_2},
\end{equation}
where we have introduced a compact notation for the 
sensitivities arriving from the previous layer:
\begin{equation}\label{eq:PQ_def}
P^{(1)}_2 = \frac{\partial a^{(1)}_2}%
{\partial W^{(1)}_{2,1}} 
= \phi'(z^{(1)}_2)\,a^{(0)},
\qquad
Q^{(1)}_2 = \frac{\partial \dot{a}^{(1)}_2}%
{\partial W^{(1)}_{2,1}},
\end{equation}
with $Q^{(1)}_2$ given by \cref{eq:pr_layer1}. The 
notation is suggestive: $P$ carries the standard chain 
rule sensitivity, while $Q$ carries the additional 
sensitivity from the time derivative --- including the 
product rule contribution from layer 1. Notice how 
$Q^{(1)}_2$ from the first product rule now feeds into 
the second product rule at layer 2 through \cref{eq:pr_layer2}. 
The product rules \emph{accumulate} as we move forward 
through the network.

\paragraph{Assembling the result.}
The full derivative of $\hat{y}'$ with respect to 
$W^{(1)}_{2,1}$ sums over all three neurons of layer 2:
\begin{equation}\label{eq:dyprime_dW1}
\frac{\partial \hat{y}'}{\partial W^{(1)}_{2,1}} = 
\sum_{k=1}^{3} W^{(3)}_{1,k}\,
\frac{\partial \dot{a}^{(2)}_k}%
{\partial W^{(1)}_{2,1}},
\end{equation}
where each term in the sum contains the product rule from 
\cref{eq:pr_layer2}, which itself carries the accumulated 
product rule from \cref{eq:pr_layer1}.

\subsubsection{The General Pattern}

The first hidden layer example reveals a pattern that 
extends to networks of any depth. We now state it 
precisely as a set of recursive relations.

For any parameter $\theta$ at layer $\ell$, the gradient 
of the total loss requires $\partial \hat{y}/\partial 
\theta$ and $\partial \hat{y}'/\partial \theta$. To 
compute these, we define two sensitivity quantities at 
each neuron $j$ of layer $m$:
\begin{equation}\label{eq:gen_PQ}
P^{(m)}_j = \frac{\partial a^{(m)}_j}{\partial \theta}, 
\qquad 
Q^{(m)}_j = \frac{\partial \dot{a}^{(m)}_j}%
{\partial \theta}.
\end{equation}

\paragraph{At the output layer ($m = L$, identity 
activation).}
The quantities we seek are:
\begin{equation}\label{eq:gen_output}
\frac{\partial \hat{y}}{\partial \theta} = 
\sum_{j} W^{(L)}_{1,j}\, P^{(L-1)}_j, 
\qquad 
\frac{\partial \hat{y}'}{\partial \theta} = 
\sum_{j} W^{(L)}_{1,j}\, Q^{(L-1)}_j.
\end{equation}
So we need $P^{(L-1)}_j$ and $Q^{(L-1)}_j$ --- the 
sensitivities at the last hidden layer.

\paragraph{At the last hidden layer ($m = L-1$).}
This layer uses a non-linear activation, so the 
sensitivities involve $\phi'$ and $\phi''$. They are 
expressed in terms of those arriving from the previous 
layer:
\begin{align}
P^{(L-1)}_j &= \phi'(z^{(L-1)}_j) 
\sum_{k} W^{(L-1)}_{j,k}\, P^{(L-2)}_k,
\label{eq:gen_P_Lm1} \\[4pt]
Q^{(L-1)}_j &= 
\phi''(z^{(L-1)}_j) 
\left(\sum_{k} W^{(L-1)}_{j,k}\, P^{(L-2)}_k\right) 
\dot{z}^{(L-1)}_j 
\;+\; 
\phi'(z^{(L-1)}_j) 
\sum_{k} W^{(L-1)}_{j,k}\, Q^{(L-2)}_k.
\label{eq:gen_Q_Lm1}
\end{align}
So we need $P^{(L-2)}_k$ and $Q^{(L-2)}_k$. The pattern 
is now clear: each layer's sensitivities depend on those 
of the previous layer.

\paragraph{General recursive relations 
($m = \ell+1, \ldots, L-1$).}
The relations above generalize to any hidden layer $m$:
\begin{align}
P^{(m)}_j &= \phi'(z^{(m)}_j) 
\sum_{k} W^{(m)}_{j,k}\, P^{(m-1)}_k,
\label{eq:gen_P_prop} \\[4pt]
Q^{(m)}_j &= 
\phi''(z^{(m)}_j) 
\left(\sum_{k} W^{(m)}_{j,k}\, P^{(m-1)}_k\right) 
\dot{z}^{(m)}_j 
\;+\; 
\phi'(z^{(m)}_j) 
\sum_{k} W^{(m)}_{j,k}\, Q^{(m-1)}_k.
\label{eq:gen_Q_prop}
\end{align}
These are applied from layer $L-1$ backward to layer 
$\ell+1$, where they meet the initialization at the 
parameter's own layer.

\paragraph{Initialization at the parameter's layer 
($m = \ell$).}
For a weight $\theta = W^{(\ell)}_{i,j}$, only 
neuron $i$ has nonzero sensitivities:
\begin{align}
P^{(\ell)}_{i} &= \phi'(z^{(\ell)}_{i})\, 
a^{(\ell-1)}_{j}, 
\label{eq:gen_P_init} \\[4pt]
Q^{(\ell)}_{i} &= 
\phi''(z^{(\ell)}_{i})\, a^{(\ell-1)}_{j}\, 
\dot{z}^{(\ell)}_{i} 
\;+\; 
\phi'(z^{(\ell)}_{i})\, \dot{a}^{(\ell-1)}_{j},
\label{eq:gen_Q_init}
\end{align}
with $P^{(\ell)}_n = Q^{(\ell)}_n = 0$ for all 
$n \neq i$. For a bias $\theta = b^{(\ell)}_{i}$, 
replace $a^{(\ell-1)}_{j}$ with $1$ and 
$\dot{a}^{(\ell-1)}_{j}$ with $0$.

\paragraph{Total gradient.}
Combining \cref{eq:gen_output} with the loss function, 
the gradient with respect to any parameter $\theta$ is:
\begin{equation}\label{eq:gen_grad_total}
\frac{\partial L}{\partial \theta} = 
2R \left( 
\frac{\partial \hat{y}}{\partial \theta} + 
\frac{\partial \hat{y}'}{\partial \theta} 
\right)_{\!t=t_c}
\!\!+ \;\lambda\, 2\bigl(\hat{y}(0)-1\bigr) 
\left(\frac{\partial \hat{y}}{\partial \theta}
\right)_{\!t=0},
\end{equation}
where the initial condition loss contributes only through 
$P$ (since $L_{IC}$ depends on $\hat{y}(0)$ but not on 
$\hat{y}'(0)$).

In practice, none of these sensitivities need to be computed 
explicitly. Frameworks such as PyTorch implement automatic 
differentiation (commonly called \emph{autograd}): during 
the forward pass, every operation is recorded in a 
computational graph along with the intermediate quantities 
--- activations, pre-activations, and their derivatives 
--- that appear in the recursive relations above. Calling 
\texttt{loss.backward()} then traverses this graph in 
reverse, combining these stored quantities to compute 
$\nabla_{\theta} L$ for all parameters simultaneously. 
Higher-order time derivatives, needed for second-order 
ODEs or PDEs, are obtained by applying autograd to the 
output of a previous autograd call, extending the 
computational graph with each application. The companion 
notebook verifies this: every gradient computed by hand 
in this section is reproduced identically by PyTorch's 
autograd engine.

With the gradient computation now established both by 
example and in full generality, a natural question 
remains: does this training cycle actually converge to the 
correct solution? The next section answers this question.


\section{Validation: Comparison with the Analytical Solution}
\label{sec:validation}

The preceding sections dissected a single iteration of the PINN 
training cycle in full numerical detail: one forward pass to 
evaluate the composite loss, one backward pass to compute 
parametric gradients, and one parameter update via gradient 
descent. At this point you might reasonably ask: does this 
actually work? If we repeat this forward--backward--update cycle 
thousands of times, will the network actually find the right 
solution?

This section answers that question affirmatively. The 1-3-3-1 
MLP is trained to convergence and its prediction $\hat{y}(t)$ 
is compared against the known analytical solution 
$y(t) = e^{-t}$ over the continuous domain $t \in [0, 1]$. The 
goal is to demonstrate that minimizing a loss function built 
exclusively from the ODE residual and the initial condition --- 
without any data sampled from the true solution --- is 
sufficient to recover the correct physical behavior.

\subsection{Training Configuration}
\label{subsec:training_config}

To transition from the single-step didactic example to a fully 
trained model, we retain the same 1-3-3-1 MLP architecture and 
the $\tanh$ activation function from 
\Cref{sec:explanatory_example}. However, the manually chosen 
parameters from 
\Cref{tab:params_l1,tab:params_l2,tab:params_l3} --- which 
served the specific purpose of enabling hand-verifiable 
calculations --- are now replaced by a \textbf{standard random 
initialization}, as would be used in any practical application. 
Specifically, we use PyTorch's default initialization scheme 
(Kaiming uniform) with a fixed random seed for full 
reproducibility. The training configuration is:

\begin{itemize}
    \item \textbf{Collocation points:} $N_c = 30$ uniformly 
    spaced points in $t \in [0, 1]$.
    
    \item \textbf{Optimizer:} Adam \cite{kingma2014adam} with 
    learning rate $\eta = 10^{-3}$. This adaptive optimizer 
    adjusts step sizes per-parameter, leading to faster and 
    more stable convergence than the basic gradient descent 
    used in the single-step example 
    (\Cref{sec:backpropagation}).

    \item \textbf{Training epochs:} $15{,}000$ complete 
    iterations of the forward--backward--update cycle.

    \item \textbf{Loss weighting:} $\lambda = 10$, consistent 
    with the value used in the single-step example.
\end{itemize}

The composite loss function remains the same as in 
\cref{eq:total_loss}:
\begin{equation}\label{eq:total_loss_val}
L = L_R + \lambda \cdot L_{IC}
  = \frac{1}{N_c}\sum_{i=1}^{N_c}
    \bigl[\hat{y}'(t_i) + \hat{y}(t_i)\bigr]^2
  + \lambda \bigl(\hat{y}(0) - 1\bigr)^2,
\end{equation}
where the temporal derivative $\hat{y}'(t_i)$ at each 
collocation point is computed via PyTorch's automatic 
differentiation engine with respect to the network 
input~$t$. Note that this is the same derivative that was 
computed analytically via forward-mode differentiation in 
\Cref{sec:explanatory_example}; the two approaches yield 
identical numerical results, as the reader can verify in the 
companion notebook.

\subsection{Evaluation Protocol and Error Metrics}
\label{subsec:eval_metrics}

After training, the network is evaluated on a fine grid of 
$N_{\text{eval}} = 500$ uniformly spaced points in $[0, 1]$. 
This grid is deliberately denser than the $N_c = 30$ 
collocation points used during training, so that the comparison 
also tests the network's ability to \textbf{generalize} to 
unseen locations within the domain. At each evaluation 
point~$t_i$, the PINN prediction $\hat{y}(t_i)$ is compared 
with the analytical solution $y(t_i) = e^{-t_i}$ using two 
standard error metrics.

\paragraph{Mean Squared Error (MSE).}
This metric quantifies the average squared deviation:
\begin{equation}\label{eq:mse}
\text{MSE} = \frac{1}{N_{\text{eval}}}
             \sum_{i=1}^{N_{\text{eval}}}
             \bigl(\hat{y}(t_i) - y(t_i)\bigr)^2.
\end{equation}

\paragraph{Relative $L^2$ Error.}
This dimensionless metric normalizes the error by the magnitude 
of the true solution:
\begin{equation}\label{eq:rel_L2}
\epsilon_{L^2}
  = \frac{\left\| \hat{y} - y \right\|_2}
         {\left\| y \right\|_2}
  = \frac{\displaystyle\sqrt{\sum_{i=1}^{N_{\text{eval}}}
    \bigl(\hat{y}(t_i) - y(t_i)\bigr)^2}}
   {\displaystyle\sqrt{\sum_{i=1}^{N_{\text{eval}}} 
    y(t_i)^2}}.
\end{equation}

\subsection{Results}
\label{subsec:results}

\paragraph{Training Convergence.}
\Cref{tab:training_evolution} summarizes the evolution of the 
loss components during training. The total loss decreases by 
approximately 6.3 orders of magnitude, from 
$1.213 \times 10^{1}$ at the first epoch to 
$6.273 \times 10^{-6}$ at epoch $15{,}000$. A clear two-phase 
dynamic is visible. In the first phase, the heavily weighted 
initial condition loss $L_{IC}$ dominates the optimization: it 
drops from $1.212 \times 10^{0}$ to below $10^{-5}$ by epoch 
$2{,}090$, meaning the network learns the correct value at 
$t = 0$ within the first $\sim\!14\%$ of training. In the 
second phase, with $L_{IC}$ effectively satisfied, the 
optimizer dedicates its full effort to reducing $L_R$ across 
the domain. By the end of training, the residual loss accounts 
for nearly all of the total loss 
($L_R = 6.210 \times 10^{-6}$ versus 
$L_{IC} = 6.341 \times 10^{-9}$, a ratio of approximately 
979:1).

\begin{table}[htbp]
\centering
\footnotesize
\caption{Evolution of the loss components during training.}
\label{tab:training_evolution}
\begin{tabular}{@{}rccc@{}}
\toprule
\textbf{Epoch} & \textbf{Total Loss $L$} & $\boldsymbol{L_R}$
               & $\boldsymbol{L_{IC}}$ \\
\midrule
     1  & $1.213 \times 10^{1}$  & $1.096 \times 10^{-2}$ 
        & $1.212 \times 10^{0}$   \\
3\,000  & $7.954 \times 10^{-3}$ & $7.911 \times 10^{-3}$ 
        & $4.251 \times 10^{-6}$  \\
6\,000  & $3.069 \times 10^{-4}$ & $3.067 \times 10^{-4}$ 
        & $2.845 \times 10^{-8}$  \\
9\,000  & $3.609 \times 10^{-5}$ & $3.609 \times 10^{-5}$ 
        & $5.914 \times 10^{-10}$ \\
12\,000 & $1.432 \times 10^{-5}$ & $1.432 \times 10^{-5}$ 
        & $3.139 \times 10^{-11}$ \\
15\,000 & $6.273 \times 10^{-6}$ & $6.210 \times 10^{-6}$ 
        & $6.341 \times 10^{-9}$  \\
\bottomrule
\end{tabular}
\end{table}

\Cref{fig:loss_history} visualizes this convergence on a 
logarithmic scale. The dominance of $L_R$ in the final training 
phase is clearly visible: after approximately epoch $3{,}000$, 
the total loss curve and the $L_R$ curve become virtually 
indistinguishable, confirming that $L_{IC}$ has been fully 
satisfied. It is worth noting that $L_{IC}$ exhibits a 
non-monotonic behavior at late epochs: after reaching a minimum 
of $3.139 \times 10^{-11}$ at epoch $12{,}000$, it increases 
slightly to $6.341 \times 10^{-9}$ at epoch $15{,}000$. This 
fluctuation is a known consequence of multi-objective 
optimization in PINNs \cite{wang2021understanding}: as the 
optimizer continues to reduce $L_R$, it may marginally 
compromise $L_{IC}$, since the two loss components share the 
same set of parameters. In this case, the fluctuation is 
entirely inconsequential --- $L_{IC}$ remains nine orders of 
magnitude below its initial value, and its contribution to the 
total loss is negligible.

\begin{figure}[htbp]
\centering
\includegraphics[width=0.65\textwidth]{figures/pinn_loss_history.pdf}
\caption{Evolution of the total loss $L$ and its components 
($L_R$ and $\lambda \cdot L_{IC}$) during $15{,}000$ training 
epochs, displayed on a logarithmic scale. The total loss 
decreases by approximately 6.3 orders of magnitude. After the 
initial phase, $L_R$ dominates the optimization.}
\label{fig:loss_history}
\end{figure}

\begin{figure}[htbp]
\centering
\includegraphics[width=\textwidth]{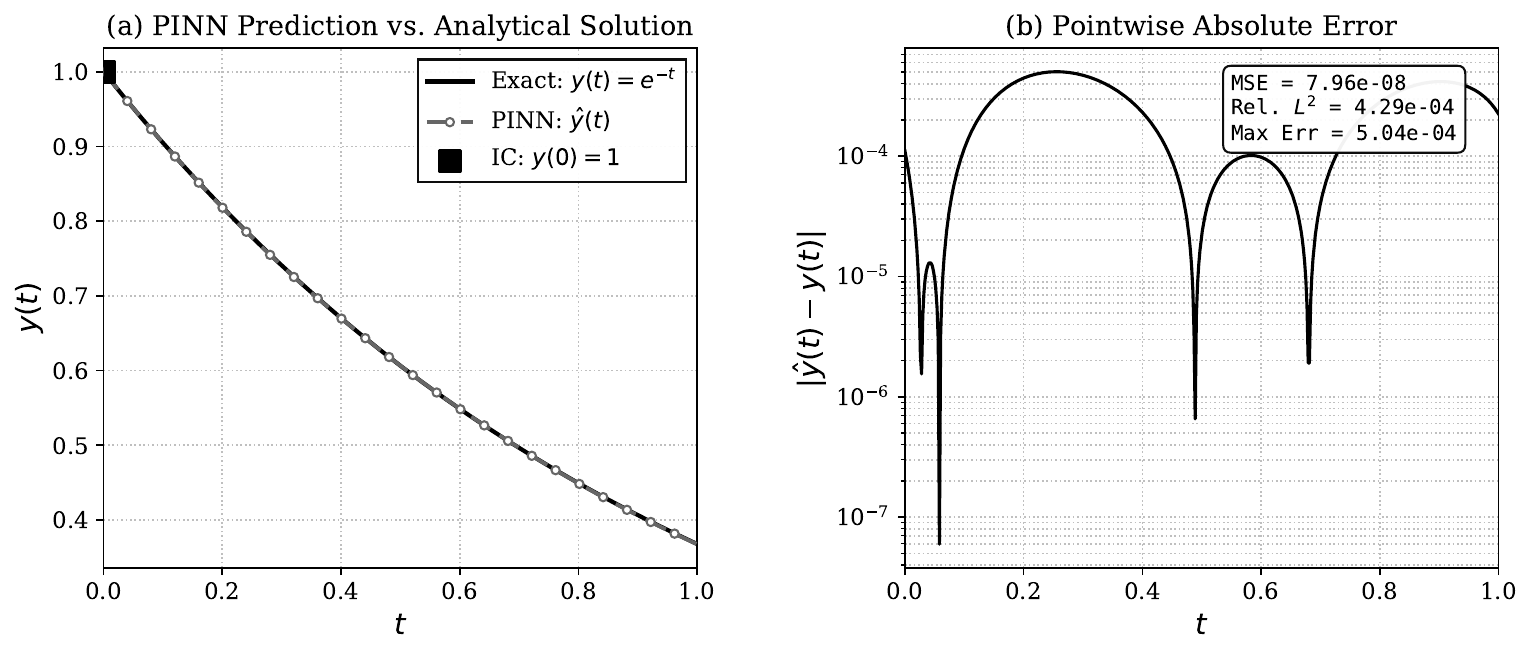}
\caption{Validation of the trained PINN against the analytical 
solution $y(t) = e^{-t}$. \textbf{(a)}~PINN prediction 
$\hat{y}(t)$ (dashed line with circle markers) overlaid on the 
exact solution (solid line); the curves are visually 
indistinguishable. \textbf{(b)}~Pointwise absolute error 
$|\hat{y}(t) - y(t)|$ on a logarithmic scale. The maximum 
error of $5.038 \times 10^{-4}$ occurs at $t \approx 0.26$.}
\label{fig:validation_plot}
\end{figure}

\paragraph{Solution Accuracy.}
\Cref{fig:validation_plot} presents the central result. The 
left panel displays the PINN prediction $\hat{y}(t)$ 
superimposed on the analytical solution $y(t) = e^{-t}$, 
evaluated on a fine grid of $N_{\text{eval}} = 500$ points. 
The two curves are visually indistinguishable over the entire 
domain.

The right panel quantifies this agreement by displaying the 
pointwise absolute error $|\hat{y}(t_i) - y(t_i)|$ on a 
logarithmic scale. The maximum error of 
$5.038 \times 10^{-4}$ occurs at $t \approx 0.26$. This 
spatial error distribution is noteworthy: the region near 
$t = 0$ benefits from the strong anchor provided by $L_{IC}$, 
while the interior of the domain relies entirely on the 
$N_c = 30$ collocation points to enforce the ODE residual. The 
error peak in the intermediate region reflects the transition 
between these two sources of constraint. The mean absolute 
error across the domain is $2.294 \times 10^{-4}$, with a 
standard deviation of $1.643 \times 10^{-4}$.

\Cref{tab:error_metrics} summarizes the global error metrics, 
and \Cref{tab:pointwise_predictions} provides the network 
predictions at five representative points.

\begin{table}[htbp]
\centering
\footnotesize
\caption{Error metrics for the trained PINN evaluated on 
$N_{\text{eval}} = 500$ points in $[0,1]$.}
\label{tab:error_metrics}
\begin{tabular}{@{}lc@{}}
\toprule
\textbf{Metric} & \textbf{Value} \\
\midrule
Mean Squared Error (MSE)       & $7.961 \times 10^{-8}$ \\
Relative $L^2$ Error           & $4.290 \times 10^{-4}$ \\
Maximum Absolute Error         & $5.038 \times 10^{-4}$ \\
\bottomrule
\end{tabular}
\end{table}

\begin{table}[htbp]
\centering
\footnotesize
\caption{PINN predictions and absolute errors at selected 
points.}
\label{tab:pointwise_predictions}
\begin{tabular}{@{}cccc@{}}
\toprule
$t$ & $\hat{y}(t)$ \textbf{(PINN)} 
    & $y(t) = e^{-t}$ \textbf{(Exact)}
    & \textbf{Absolute Error} \\
\midrule
$0.00$ & $0.999889$ & $1.000000$ & $1.109 \times 10^{-4}$ \\
$0.25$ & $0.777907$ & $0.778411$ & $5.032 \times 10^{-4}$ \\
$0.50$ & $0.605947$ & $0.605923$ & $2.402 \times 10^{-5}$ \\
$0.75$ & $0.471494$ & $0.471657$ & $1.631 \times 10^{-4}$ \\
$1.00$ & $0.367655$ & $0.367879$ & $2.244 \times 10^{-4}$ \\
\bottomrule
\end{tabular}
\end{table}

These results confirm that the 1-3-3-1 MLP, despite its 
minimal capacity of only 22 trainable parameters, achieves a 
relative $L^2$ error of $0.043\%$ in approximating the true 
solution. It bears repeating that at no point during the 
$15{,}000$ training epochs was the model given access to any 
value of the function $y(t) = e^{-t}$. The physics encoded in 
the loss function was the sole source of information guiding 
the optimization.

\subsection{Discussion and Scope of the Validation}
\label{subsec:discussion}

Let us be clear about what this result does and does not show.

The first-order linear ODE $y'(t) + y(t) = 0$ is an elementary 
problem for which classical numerical methods --- such as 
fourth-order Runge-Kutta schemes \cite{burden2015numerical} --- 
provide solutions with machine precision at a fraction of the 
computational cost required by the PINN training loop. For this 
class of well-posed, low-dimensional initial value problems, 
PINNs are \emph{not} competitive with established solvers, nor 
is it the purpose of this didactic example to suggest otherwise.

The value of this validation is, instead, \textbf{conceptual 
and pedagogical}. It demonstrates the core mechanism of the 
PINN framework: a neural network, guided solely by the 
structure of a differential equation and an initial condition, 
can discover the correct physical solution without ever being 
shown a single data point from that solution. This ``learning 
from physics alone'' paradigm is the foundational building 
block upon which all advanced PINN applications are 
constructed.

The true potential of the PINN methodology emerges in problem 
classes where classical methods face significant limitations:

\begin{itemize}
    \item \textbf{Inverse problems:} Identifying unknown 
    parameters within the governing equation --- such as 
    diffusion coefficients, reaction rates, or material 
    properties --- from sparse and potentially noisy 
    experimental observations. In an inverse PINN formulation, 
    these physical parameters are simply added to the set of 
    trainable variables $\theta$ and are learned jointly with 
    the network weights during optimization 
    \cite{raissi2019physics, karniadakis2021physics}.

    \item \textbf{Data-physics fusion:} Combining incomplete 
    or noisy experimental measurements with known physical 
    laws to reconstruct full-field solutions. A landmark 
    example is the reconstruction of hidden velocity and 
    pressure fields from flow visualizations 
    \cite{raissi2020hidden}, a task for which no classical 
    solver can be directly applied without full boundary 
    condition specification.

    \item \textbf{High-dimensional and complex-geometry 
    problems:} Partial differential equations defined on 
    irregular domains or in high-dimensional parameter spaces, 
    where mesh generation is prohibitive or where classical 
    methods suffer from the curse of dimensionality. Because 
    PINNs are mesh-free, they circumvent the need for domain 
    discretization entirely.
\end{itemize}

\noindent In each of these settings, the training mechanics 
remain \emph{exactly} those detailed in this paper: a composite 
loss function is constructed from the governing equations and 
available data, and backpropagation adjusts the network 
parameters to satisfy all imposed constraints simultaneously. 
By mastering these mechanics on the simple problem treated 
here, the reader is equipped with the foundational 
understanding necessary to approach these more challenging and 
impactful applications with confidence.


\section{Conclusion}
\label{sec:conclusion}

This paper has presented a didactic and self-contained 
introduction to the complete training cycle of a 
Physics-Informed Neural Network, applied to a first-order 
initial value problem. Existing PINN tutorials and guides 
either delegate the internal training mechanics to 
automatic differentiation libraries or, in the case of 
Almqvist \cite{almqvist2021fundamentals}, derive them 
explicitly but for a single hidden layer only. The present 
work extends this explicit treatment to a multi-layer 
network, where a richer gradient structure emerges. Through 
a fully verifiable numerical example and a subsequent 
validation against the analytical solution 
$y(t) = e^{-t}$ (\Cref{sec:validation}), we gave 
explicit treatment to four aspects of the training process 
that are not readily visible in library-based 
implementations:

\begin{enumerate}
    \item \textbf{Dual forward propagation:} During the 
    forward pass, the network output $\hat{y}$ and its 
    temporal derivative $\hat{y}'$ are propagated 
    simultaneously through the network, layer by layer, 
    using the chain rule at each stage. This dual 
    propagation is what makes the additional gradient 
    structure of the backward pass visible.

    \item \textbf{The product rule in backpropagation:} 
    For hidden layer parameters, the temporal derivative 
    $\hat{y}'$ depends on both the activations and their 
    derivatives through the same parameter, requiring the 
    product rule --- and introducing the second derivative 
    of the activation function $\phi''$ even for a 
    first-order ODE. This is the structural complication 
    that distinguishes PINN backpropagation from standard 
    neural network backpropagation.

    \item \textbf{The role of the hyperparameter 
    $\lambda$:} The validation in \Cref{sec:validation} 
    provided direct empirical evidence of a two-phase 
    training dynamic: the $\lambda = 10$ weighting drove 
    $L_{IC}$ below $10^{-5}$ within the first $2{,}090$ 
    epochs, after which the optimizer redirected its full 
    effort toward minimizing the ODE residual across the 
    domain.

    \item \textbf{General gradient formulas:} The product 
    rule structure observed in the worked examples was 
    generalized into recursive relations for two 
    sensitivities --- $P^{(m)}_j$ for the activation and 
    $Q^{(m)}_j$ for its time derivative --- that propagate 
    from the parameter's layer through successive layers 
    to the output, accumulating product rules and fan-out 
    summations at each stage. These relations extend the 
    gradient computation to networks of arbitrary depth 
    without requiring case-by-case derivations, and they 
    correspond directly to the intermediate quantities 
    that automatic differentiation engines store and 
    combine during \texttt{loss.backward()}.
\end{enumerate}

The validation confirmed that the minimization of the 
physics-informed loss --- without access to any data from 
the true solution --- is sufficient to recover the correct 
physical behavior, achieving a relative $L^2$ error of 
$4.290 \times 10^{-4}$ ($0.043\%$) with a minimal 
22-parameter network. While PINNs are not competitive with 
classical solvers for this class of elementary problems, 
the exercise establishes the conceptual foundation 
necessary for tackling the problem classes where the 
framework offers unique advantages, namely inverse 
problems, data-physics fusion, and high-dimensional PDEs 
(\Cref{subsec:discussion}).

To bridge the gap between mathematical derivation and 
practical implementation, this work is accompanied by an 
interactive Jupyter Notebook that reproduces every manual 
calculation presented in the text using PyTorch's 
automatic differentiation engine, and independently 
replicates the full training and validation pipeline. The 
exact agreement between the hand-derived gradients and 
PyTorch's autograd output serves as mutual validation: the 
theory verifies the code, and the code verifies the 
theory.

As natural extensions of this didactic framework, we 
envision:

\begin{itemize}
    \item \textbf{Higher-order and multi-dimensional 
    problems:} Generalizing the didactic analysis to 
    second-order ODEs and to PDEs in two spatial 
    dimensions, where the interplay between multiple 
    differential operators and boundary conditions 
    introduces additional complexity in the loss function 
    and gradient computation.

    \item \textbf{Training dynamics and optimization:} 
    Investigating the impact of alternative activation 
    functions (such as GELU and Swish) on the gradient 
    flow through the physics-informed loss, and comparing 
    the convergence behavior of first-order optimizers 
    (Adam \cite{kingma2014adam}) with second-order methods 
    (L-BFGS \cite{liu1989limited}), which have shown 
    strong performance in PINN training 
    \cite{raissi2019physics}.

    \item \textbf{Inverse problems:} Extending the 
    framework to identify unknown physical parameters 
    from experimental or observational data, where the 
    unknown parameters are added to the trainable set 
    $\theta$ and learned jointly with the network 
    weights. Benchmarking against classical estimation 
    methods on well-characterized test problems would 
    clarify both the advantages and limitations of the 
    physics-informed approach.
\end{itemize}

\noindent In each case, the foundational training 
mechanics and the general gradient formulas presented here 
carry over directly. We hope that this work lowers the 
barrier to entry into Scientific Machine Learning and 
encourages the adoption of physics-informed approaches in 
both research and education.


\section*{Code and Data Availability}

This work does not involve experimental or observational 
datasets. All results are generated computationally from the 
governing differential equation and its initial condition. To 
guarantee full reproducibility, the following materials are 
provided:

\begin{enumerate}
    \item \textbf{Source code:} The complete Python/PyTorch 
    implementation is provided as an interactive Jupyter 
    Notebook, structured in two parts:
    \begin{itemize}
        \item \textit{Verification of manual calculations 
        (\Cref{sec:explanatory_example,sec:backpropagation}):} 
        The notebook algorithmically reproduces every 
        hand-derived value presented in the text using 
        PyTorch's automatic differentiation engine, serving as 
        a didactic bridge between theory and implementation.

        \item \textit{Training and validation 
        (\Cref{sec:validation}):} The notebook executes the 
        full PINN training pipeline, saves the trained model 
        checkpoint and all evaluation metrics, and generates 
        the publication figures.
    \end{itemize}

    \item \textbf{Computational artifacts:} The trained model 
    checkpoint (\url{pinn_trained_model.pt}), evaluation 
    metrics and predictions (\url{pinn_results.json}), and 
    epoch-by-epoch loss history 
    (\url{pinn_loss_history.npz}) are provided as 
    supplementary files. All numerical values reported in 
    \Cref{sec:validation} are read directly from these saved 
    files, ensuring that the manuscript numbers remain identical 
    regardless of hardware or re-execution.
\end{enumerate}

\noindent The complete codebase, including all computational 
artifacts, is publicly available at 
\url{https://github.com/Tahimi/PINN-Didactic-Training-Cycle} 
and archived with a permanent DOI via Zenodo 
\cite{PINN_notebook_2026}.


\section*{Acknowledgments}

During the development of this work, the authors used 
AI-based language models (Claude, Anthropic, 2024--2025) 
as assistive tools for writing, code development, and 
computational verification. All conceptual and scientific 
decisions --- including the identification of the 
pedagogical gap, the design of the illustrative example, 
the selection of which aspects of the training cycle to 
emphasize, and the critical discussion of the method's 
scope and limitations --- originated with the authors. 
All mathematical derivations, numerical values, and 
computational results were independently verified by the 
authors, who assume full responsibility for the integrity, 
accuracy, and originality of the submitted work.

\bibliographystyle{siamplain}
\bibliography{references}

\end{document}